\documentstyle[11pt]{article}
\begin{document}
\bibliographystyle{alpha}
\baselineskip=5.0mm

\begin{center}
{\large\bf A Counting Proof of the Graham Pollak Theorem}
\ \\
\ \\
   {\large Sundar Vishwanathan} \\
    Department of Computer Science  and Engineering \\
  Indian Institute of Technology, Bombay   \\
  Powai, Mumbai \\
  India 400076  \\
  {\tt sundar@cse.iitb.ernet.in}\\

\end{center}

\def\Box{{\rm \mbox{\,\rule[1.5pt]{.1pt}{4pt}\rule[1.5pt]{4pt}{.1pt}\hskip -4pt\rule[5.5pt]{4pt}{.1pt}\rule[1.5pt]{.1pt}{4pt}}}\,\,} 
\newtheorem{defi}{Definition}
\newenvironment{definition}{\begin{defi}\rm}{\end{defi}}
\newtheorem{claim}{Claim}

\newtheorem{theorem}{Theorem}
\newtheorem{lemma}[theorem]{Lemma}
\newtheorem{fact}[theorem]{Fact}
\newtheorem{proposition}[theorem]{Proposition}
\newtheorem{note}[theorem]{Note}
\newtheorem{corollary}[theorem]{Corollary}
\newtheorem{example}[theorem]{Example}
\newtheorem{conjecture}[theorem]{Conjecture}
\newcommand{\blackslug}{\hbox{\hskip 1pt \vrule width 4pt height 8pt
depth 1.5pt \hskip 1pt}}
\newcommand{\QED}{\quad\blackslug\lower 8.5pt\null}
\newcommand{\dg}{DAG}
\newcommand{\nullset}{\emptyset}
\newcommand{\assign}{\leftarrow}
\newenvironment{proof}{\noindent {\bf Proof:}}{\QED}
\newenvironment{sketch}{\noindent {\bf Proof (sketch):}}{\QED}
\newcommand{\nfont}{\sl}
\newcommand{\st}{\,|\;}
\newcommand{\inter}{\cap }
\newcommand{\union}{\cup }
\newcommand{\edge}{\relbar\joinrel\joinrel\relbar}
\newcommand{\qed}{\quad\blackslug\lower 8.5pt\null}
\newcommand{\nospaceqed}{\blackslug\lower 8.5pt\null}
\newcommand{\floor}[1]{\left\lfloor #1 \right\rfloor}
\newcommand{\ceil}[1]{\left\lceil  #1 \right\rceil}
\newcommand{\suchthat}{\colon\,}
\newcommand{\problem}[3]{\medskip\noindent {\bf #1}:\par\noindent
${\underline {\rm Input:}}$ #2\par\noindent ${\underline {\rm Question:}}$ #3} 

\newcommand{\does}[2]{{#1}: & {#2}.\\ }
\newcommand{\sets}[3]{{#1} & {#2} \leftarrow {#3}.\\ }
\newcommand{\sends}[3]{{#1} \rightarrow {#2}: & {#3}.\\ }
\newcommand{\domain}{\mbox{\it Dom}}
\newcommand{\range}{\mbox{\it Range}}

\newcommand{\half}{\mbox{$\frac{1}{2}$}}       
\newcommand{\inverse}[1]{\frac{1}{{#1}}}       
\newcommand{\inv}[1]{{#1}^{-1}}                
\newcommand{\prob}[1]{\mbox{Pr}\left[{#1}\right]}     
\newcommand{\mod}{\mbox{~mod~}}                
\newcommand{\ord}{\mbox{~ord~}}                
\newcommand{\midbox}{\makebox[3pt]{$\mid$}}
\newcommand{\abs}[1]{
  \ifmmode {{\midbox {#1}\midbox}}
  \else \mbox{\midbox {#1}\midbox} \fi }
\newcommand{\dblmid}{\makebox[4pt]{$\mid\-\mid$}}
\newcommand{\norm}[1]{
  \ifmmode \mbox{\dblmid ${#1}$\dblmid}
  \else \mbox{\dblmid {#1}\dblmid} \fi }
\newcommand{\lpar}{\mbox{$<$}}                 
\newcommand{\rpar}{\mbox{$>$}}                 
\newcommand{\set}[1]{\ifmmode \mbox{\{{$#1$}\}}%
  \else \mbox{\{{#1}\}} \fi}                 
\newcommand{\comment}[1]{\ifmmode {\{ {\rm {#1}}\}}%
  \else {\{ {#1}\}} \fi}                      
\newcommand{\anglebrack}[1]{\ifmmode \mbox{\lpar {$#1$}\rpar}%
   \else \mbox{\lpar {#1}\rpar} \fi}           

\newenvironment{proof2}{\addvspace{\bigskipamount}\noindent{\em Proof. }}%
{\ \ $\Box$\par\addvspace{\bigskipamount}}
\def\protocol#1#2{\medbreak\noindent{\bf Protocol #1.} \par #2\par\medbreak}
\newenvironment{sketch2}{\addvspace{\bigskipamount}\noindent{\em Proof (sketch). }}%
{\ \ $\Box$\par\addvspace{\bigskipamount}}


\newcommand{\cS}{{\cal S}}
\newcommand{\nbhd}{{\hbox{\it Neighbourhood}}}

\newcommand{\ca}{{\cal A}}
\newcommand{\asp}{\hspace*{.2in}}
\newcommand{\tisp}{\hspace*{.2in}}


\newcommand{\nbp}{{\Gamma^+}}
\newcommand{\nbn}{{\Gamma^-}}

\begin{abstract}
We give a counting based proof of the Graham-Pollak theorem.
\end{abstract}

\section{Introduction}
A spectacular application of linear algebra to prove a  combinatorial
statement is the Graham-Pollak theorem \cite{GP}.
The theorem states that the
edge set of the  complete
graph $K_n$ cannot be written as the disjoint union of
$n-2$ complete bipartite graphs. The original proof used Sylvester's
law of inertia. See \cite{BF,P,Tv} for other  short proofs.
These proofs seem to  use linear algebra inherently. Combinatorialists
have often commented that a combinatorial proof for the theorem is not known.
See, for example comments regarding the problem in \cite{VW,AZ}. 

In this note we show that the linear algebra  proof of
\cite{Tv}, for instance,  can be  explained combinatorially.
The first observation is that the linear algebra part of the proof
can be replaced by a pigeon-hole argument. The other  (minor)
observation in this note is that with this in place,
one can explain the slick calculations
in the usual linear algebra proofs of this theorem
with a  longer explicit bijective argument.

\section{The Proof}
\begin{theorem}[Graham-Pollak]
Suppose that $K_n$ is obtained  as the edge-disjoint union
of $m$ complete bipartite graphs.
Then $m \geq n-1$.
\end{theorem}
{\bf Proof:}
For a contradiction, consider a covering of $K_n$ by 
complete bipartite graphs  $(L_i,R_i): 1\leq i \leq n-2$.
The vertex set is identified with $[n]$.

Consider a labeling of the  $n$ vertices, 
$\sigma: [n]\rightarrow[k]$ where $k> n^n$.
We associate a pattern  with $\sigma$,  as an $n-1$ tuple, where
the $i$th entry of  the tuple, for $i < n-1$,
is given by the sum of the $\sigma$
values of the vertices in $L_i$; that is the $i$th entry 
is $\Sigma_{j\in L_i} \sigma(j)$.
The $n-1$th entry is the sum of the $\sigma$ values of all vertices.
The number of possible patterns is at most $(kn)^{n-1}.$
The total number of labelings is $k^n$. Hence, since $k$ is large
enough, there are
two distinct labelings with the same pattern, say $\sigma_1$ and $\sigma_2$.
Define $\tau=\sigma_1-\sigma_2$. Note that for each $1\leq i\leq n-2$,
$\Sigma_{j\in L_i} \tau(j) =0.$ Also $\Sigma_{j=1}^n \tau(j)= 0$ and
$\tau$ is non-zero on at least one vertex. 

Consider the following equality:
\[ ( \Sigma_{j=1}^n \tau(j))^2=   \Sigma_{j=1}^n (\tau(j))^2 +
 2\Sigma_{i<j} \tau(i)\tau(j)  \]
The left hand side is zero.
The first term in the right hand side is non-zero. For a contradiction
we will show that the second term is zero.
Because we have a disjoint cover of $K_n$,
 \[ \Sigma_{i<j} \tau(i)\tau(j) = 
 \Sigma_{i=1}^n ( \Sigma_{j\in L_i} \tau(j)) (\Sigma_{k\in R_i} \tau(k)).\]
But the right hand side is zero since for each $i$, $\Sigma_{j\in L_i} \tau(j) =0$.
\qed

\section{Explaining the Calculations}
The last part of the proof; defining $\tau$ and  the calculations following
it seem rather mysterious. Especially the use of  the fact that if  the sum
of squares is zero then each term should be zero.
We give an explicit bijection to explain this phenomenon. We continue
the proof after defining $\sigma$. 

We consider two graphs $H$ and $H'$ defined below. Both have 
the same vertex set $W$ which is partitioned into $2n$ non-empty parts: 
$W= V_1\cup V_1' \cup V_2 \cup V_2'\cdots V_n\cup V_n'$.
We will require $\Sigma_{i=1}^n |V_i|= \Sigma_{i=1}^n |V_i'|=N$ (say.)
$|V_i|$ could  be any positive integer.
The edge set of $H$ is as follows:
$\{ uv: u\in V_i, v\in V_j, i\not=j \} \cup
\{ uv: u\in V_i', v\in V_j', i\not=j \}$. Note that $H$ is the union of two
disjoint complete $n$-partite graphs, each on $N$ vertices.
The edge set of $H'$ is as follows:
$\{ uv: u\in V_i, v\in V_j', i\not=j \}$. Hence $H'$ is a bipartite graph.
The key observation is this. 
\begin{lemma}
If the number of edges in $H$ and $H'$ are the same then 
$|V_i| = |V_i'|$ for all $ 1\leq i \leq n$.
\end{lemma}
{\bf Proof.}  (Sketch.)
We first claim that one can prove by an explicit bijection
the fact that the number of edges of
$K_p \cup K_q$ plus $\min\{p,q\}$ is at least the number of edges
in $K_{pq}$.
They are equal iff $p=q$ or $q=p+1$, assuming $q$ is larger.
To see the bijection, assume $p < q$, consider two sets of vertices
of size $p$ (say $P$) and $q$ (say $Q$).
Map an edge $(i,j)$ (for $i\in P, j\in Q, i\not=j$)
of $K_{pq}$ to an edge in $K_p \cup K_q$ as follows.
If $i< j\leq p$ then to
$ij$ in $K_p$. Otherwise to $ij$ in $K_q$. Note that this is one to one.
Also note that the edges $ij$ with both $i,j$ greater than $p$ are not in 
the range of the map. 
Hence the number of edges in $K_p \cup K_q$ plus 
$\frac{p+q}{2}$ is at least the number of edges in $K_{pq}$. 
And they are equal only if $p=q$.

If we take the complement of each component of $H$ with respect
to the complete graph $K_N$, then we get disjoint cliques of size
$p_i$ and $q_i$ for each $i$. If we take the complement of $H'$,
with respect to $K_{N,N}$, then we get disjoint copies  of $K_{p_i,q_i}$s
for each $i$.  We will use these complements in the next paragraph.

If the  number edges in $H$ and $H'$ are both $m$, then the number of edges
in the complement of $H$ is $2{N \choose 2}-m$ and that in the
complement of $H'$ is $N^2-m$. That is the former plus $N$ equals
the latter.  However, by the argument in the first  paragraph,  if 
$p_i \not= q_i$ for any one $i$, then by adding the inequality for
the  number of edges for each $i$, we infer that the number of edges 
in the complement of $H$ plus $N$ is strictly
larger than the latter, a contradiction.
\qed

Consider a disjoint cover of the edges of $K_n$  by bipartite
graphs $(A_j, B_j): 1\leq j\leq k$. Here $k$ is not restricted.    
We define $H, H'$ as a disjoint union of $n$ sets of vertices,
the $i$th (called $V_i$ and $V_i'$ in $H$ and $H'$ respectively)
corresponding to some number of copies of the $i$th vertex of $K_n$.
Now, suppose that for each $1\leq j\leq k$, 
$\Sigma_{p\in A_j} |V_p| =\Sigma_{p\in A_j} |V_p'|$ then  we claim that
the number of edges  in $H$ and $H'$ are the same.
Indeed, we can cover the edges of $H$ by the following bipartite graphs;
two for each $1\leq j\leq k$:
\[(\cup_{p \in A_j} V_p, \cup_{q \in B_j} V_q) \]  
\[(\cup_{p \in A_j} V_p', \cup_{q \in B_j} V_q') \]  
Similarly we can cover the edges of $H'$ by the following bipartite
graphs, two for each $j$:
\[(\cup_{p \in A_j} V_p, \cup_{q \in B_j} V_q') \]  
\[(\cup_{p \in A_j} V_p', \cup_{q \in B_j} V_q). \]  

It can be seen that, for each $j$,
the total number of edges in the top two bipartite graphs
covering edges in $H$, 
is equal to the total number of edges in the bottom two bipartite graphs.
This implies that the number of edges in $H$ and $H'$ are
equal. By the lemma this implies that $|V_p| = |V_p'|$, for every $p$.

To use the above,  to yield a contradiction for the Graham Pollak result,
we use $\sigma_1$ and $\sigma_2$ constructed earlier, to define these two 
graphs, with $|V_p| \not= |V_p'|$ for at least one $p$.
\qed

The main observation resulting in this note is that the following fact, which
follows from linear algebra, can also be proved by a pigeon hole argument. 
\begin{lemma}
Let $A$ be an $m\times n$ integer matrix, with $m < n$.
Then, there are positive integer vectors $x_1$ and
$x_2\not=x_1$  such that $Ax_1=Ax_2$.
\end{lemma}
Essentially, if the domain is restricted to a large enough (finite) set,
then the range can be made smaller than the domain and hence  two domain
elements map to the same point  in the range.  
Other theorems where the linear algebra part used is the lemma above
can be proved using the pigeon hole principle. 

\section{On Proofs.} What constitutes a combinatorial proof and what is
a linear algebra (or topological) proof are  questions faced by 
mathematicians for some time now.  It is quite conceivable that different
people have different notions! In the  proof in this note  we do not use 
the notion of a field, which seems to be necessary for linear algebra
based proofs. The  Steinitz exchange  lemma
and Gaussian Elimination, both use the concept of inverses in the field.

The proof in this note however  constructs intermediate structures of 
large size. But implicitly so do the linear algebra proofs. Even if
we work over the rationals, for the intermediate values, the
numerator and denominator may be as
large as the determinant of an $n\times n$ $0-1$ matrix each of which
can be as large as $n^{O(n)}$. We insist that these have to be written in
unary.  For a graph on $n$ vertices we either construct graphs on $n^n$ vertices
or use labels of size $n^n$ (again assuming labels are written in unary.)
Typically, one sees this phenomenon of using large numbers for the proof
in Ramsey Theory and not so much in other areas of Combinatorics. 
Call a proof an {\em effective} combinatorial proof if the size of the
proof (assuming that intermediate labels are written in unary) 
is polynomial in the input size. (Here the input is  an explicit description
of $K_n$.) Finding an
effective combinatorial proof for the Graham Pollak theorem is 
a nice open problem. This method will give us an effective proof for a 
worse bound, by restricting $k$ to be a polynomial
in $n$.  A bound of $n/2$ by an effective combinatorial proof which does  not 
use counting would be a nice first question to solve.
Note that a bound of $n/2$ can be proved using linear algebra over 
${\cal F}_2$, but this uses {\em counting with parity} and/or exact
counting with large numbers. 


%
%
\end{document}